\documentclass{amsart}

\usepackage{amssymb}
\usepackage{amsthm}
\usepackage{amsmath, tikz}
\usepackage{amsbsy}
\usepackage[all]{xy}
\usepackage{bm}
\usepackage{hyperref}

\begin{document}

\newtheorem*{theorem}{Theorem}
\newtheorem*{conjecture}{Conjecture}
\newtheorem{lemma}{Lemma}

\title[]{Fundamental component enhancement via adaptive nonlinear activation functions}
\thanks{S.S. was partly supported by the NSF (DMS-2123224) and the Alfred P. Sloan Foundation.}

\author[]{Stefan Steinerberger}
\address{Department of Mathematics, University of Washington, Seattle, WA 98195, USA}
\email{steinerb@uw.edu}

\author[]{Hau-tieng Wu}
\address{Department of Mathematics and Department of Statistical Science, Duke University, Durham, NC, USA}
\email{hauwu@math.duke.edu}

\begin{abstract}
In many real world oscillatory signals, the fundamental component of a signal $f(t)$ might be weak or does not exist. This makes it difficult to estimate the instantaneous frequency of the signal. Traditionally, researchers apply the rectification trick,  working with $|f(t)|$ or $\mbox{ReLu}(f(t))$ instead, to enhance the fundamental component. This raises an interesting question: what type of nonlinear function $g:\mathbb{R} \rightarrow \mathbb{R}$ has the property that
$g(f(t))$ has a more pronounced fundamental frequency? $g(t) = |t|$ and $g(t) = \mbox{ReLu}(t)$ seem to work well in practice; we propose a variant of 
$g(t) = 1/(1-|t|)$ and provide a theoretical guarantee. Several simulated signals and real signals are analyzed to demonstrate the performance of the proposed solution.
\end{abstract}

\maketitle

\section{Introduction}

\subsection{The problem.} We start by describing the problem in the simplest possible setting: suppose we are given a $1-$periodic signal
\begin{equation} \label{Simple model harmonic 1}
f(t) = \sum_{k=0}^{\infty} a_k \cos{(2\pi kt)} + b_k \sin{(2\pi kt)}
\end{equation}
and define the support of $S$ to be
$ S = \left\{k \in \mathbb{N}: a_k^2 + b_k^2 > 0 \right\}$.  Our goal is to recover $\gcd(S)$, the greatest common divisor of all the elements in $S$.
If we have recorded two or more periods of the signal, then the $\gcd$ will be bigger than 1. In that case, we call $\gcd(S)$ the fundamental
frequency. Mathematically speaking, there is no difficulty: compute the Fourier transform of $f$ and check. 

 \begin{center}
 \begin{figure}[h!]
 \includegraphics[width =0.8\textwidth]{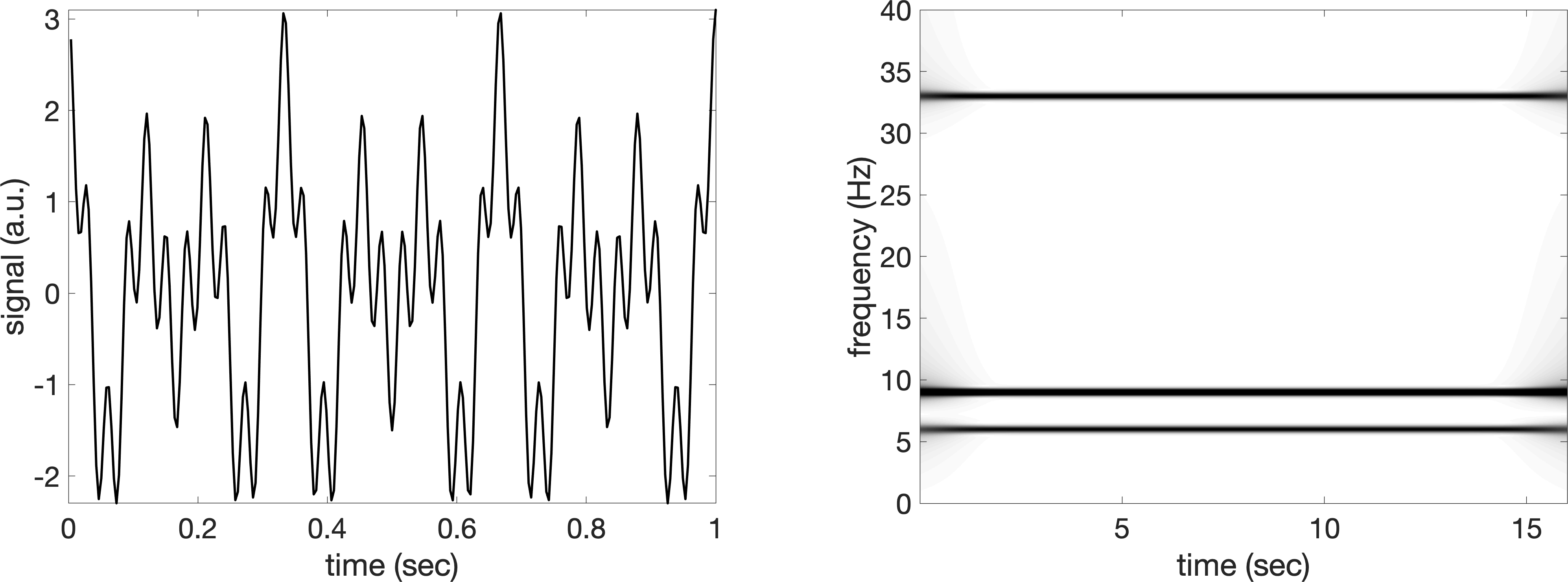}
  \vspace{-5pt}
 \caption{Left: $f(t) = 0.8 \cos{(2\pi 6t)} + 1.4\cos{(2\pi 9t)} + 0.9\cos{(2\pi 33t)}$ has $\gcd(S) = 3$.  Right: the spectrogram shows three distinct frequencies but not the fundamental frequency 3n. \label{Figure1SimpleExample}}
 \end{figure}
 \end{center}
 
There are several reasons why this seemingly trivial problem is difficult in practice.
\begin{enumerate}
\item In practice, all Fourier coefficients of $f$ will be distinct from 0; some form of thresholding to decide significance is required.
\item The actual fundamental frequency may not be constant, it may change from time to time, or it may not exist at all \cite{Colominas2021,Hou2016,Torresani2019,Wu2013,Yang2019}.
\item  We only have $f$ sampled at finitely many points; moreover, in practice one should always expect some form of noise. 
\end{enumerate}

\subsection{Motivation from medicine.} 
\begin{center}
 \begin{figure}[h!]
 \includegraphics[width =\textwidth]{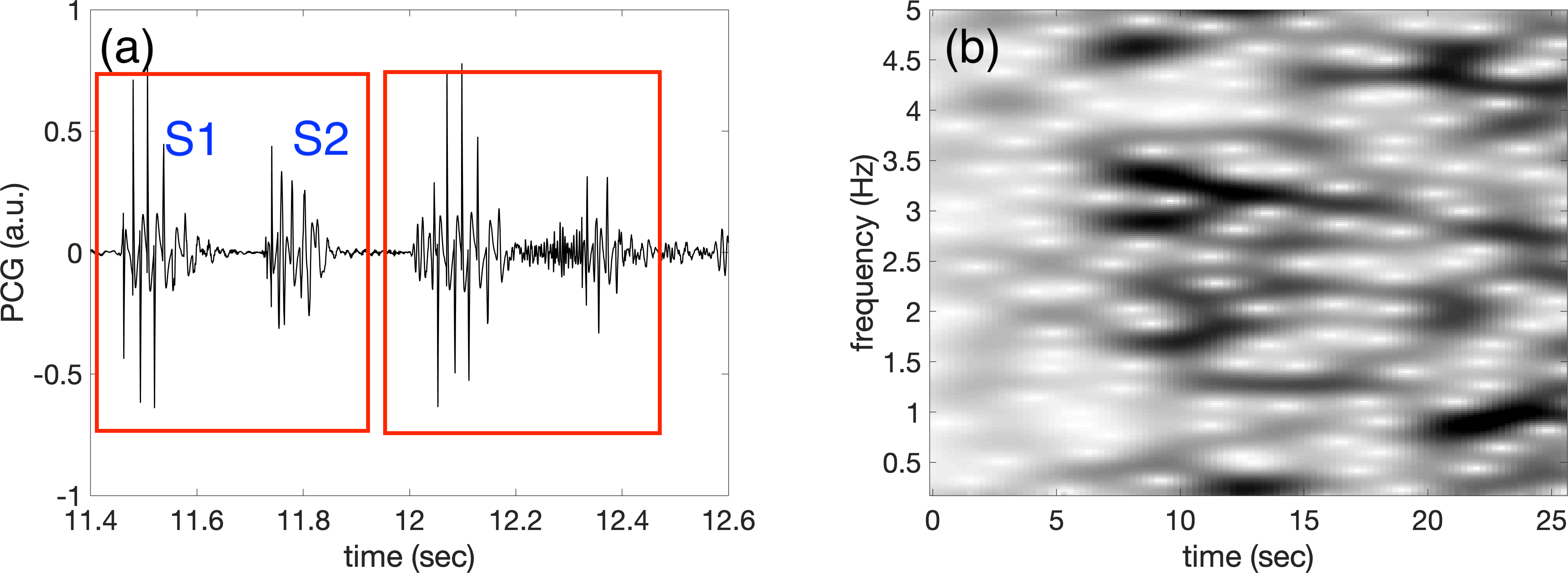}
  \vspace{-10pt}
 \caption{Phonocardiogram (a) and Spectrogram (b): no clear fundamental frequency is visible. The heart rate of this subject is about 2Hz; that is, there are two cardiac cycles per second.}
 \label{fig:pcg}
 \end{figure}
 \end{center}

We describe an explicit problem coming from medical signal processing which was partially the motivation for our paper. Our example can be seen in Figure \ref{fig:pcg}(a) and comes from a phonocardiogram (PCG) signal \cite{PCGreview2017}: this is a sound signal created by the vibrations created by the closure of the heart valves.  One fundamental component of the signal, one heartbeat, can be seen within each of the two red boxes: it is comprised of two ingredients, S1 and S2 (indicated by blue errors): S1 is due to the atrioventricular valves closing at the beginning of systole and S2 is the consequence of the aortic and pulmonary valves closing at the end of systole. 

Due to the heart rate variability, the periods between two consecutive cycles are not fixed: an important but challenging problem is to estimate the heart rate at a local point in time, this amounts to estimating the time-varying period, or the instantaneous frequency. 
Naturally, one would be inclined to apply time-frequency (TF) analysis \cite{Flandrinbook}. This signal oscillates roughly twice per second, so the frequency is around 2 Hz, and we would expect to see a dominant curve around 2 Hz that represents the fundamental component of the signal. However, we cannot see anything concrete around 2 Hz from the standard spectrogram shown in Figure \ref{fig:pcg}(b) and, in this sense, the fundamental component of the signal does not exist.
In short, we need a different solution if we want to estimate the time-varying heart rate.

%

%
%
%

\subsection{The rectification trick.} 
\begin{center}
 \begin{figure}[h!]
\includegraphics[width =0.95\textwidth]{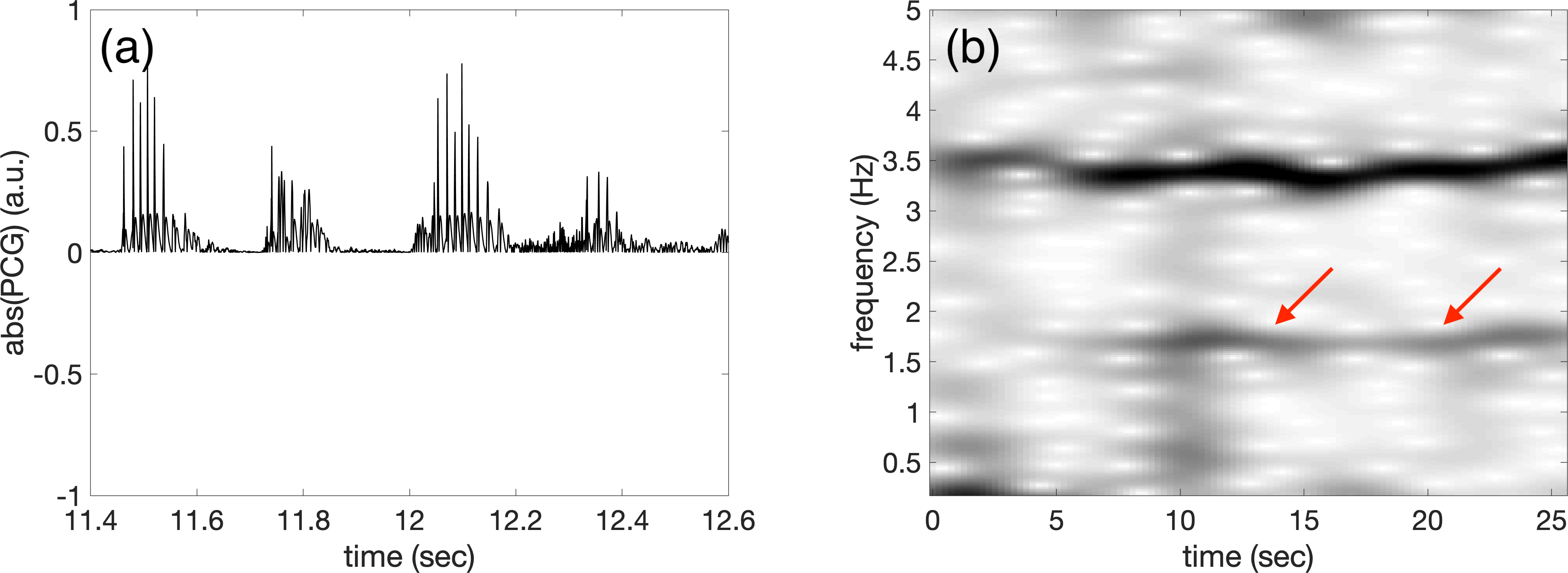}
  \vspace{-10pt}
 \caption{Rectified Phonocardiogram (a) and its Spectrogram (b): a fundamental frequency around 2 Hz starts emerging, indicated by red arrows.}
 \label{fig:rect}
 \end{figure}
 \end{center}

There is a surprisingly simple and widely applied solution; that is, take the `rectification' of the signal before running any sort of time-frequency analysis. Mathematically speaking, take the signal $\left|f(t)\right|$ where $f(t)$ is the signal of interest.
We refer to Figure \ref{fig:rect} where the signal from Figure \ref{fig:pcg} has been analyzed after rectification: we analyze $\left|f(t)\right|$ instead of $f(t)$.

  \vspace{-20pt}

   \vspace{-20pt}

Here, rectification leads to a dramatically different and much more informative spectrogram in which a clearly defined fundamental frequency exists (indicated by the red arrows). 
This simple solution has been widely applied in practice if we want to estimate the fundamental frequency, for example, in the problem of extracting fetal electrocardiogram (ECG) from the trans-abdominal maternal ECG \cite{SuMiller2019} or extracting f-wave from the ECG signal of a patient with atrial fibrillation \cite{Malik2017}, which is a far from complete list.
To the best of our knowledge, there is no theoretical argument explaining why the rectification trick works to enhance the fundamental component. We believe this to be of substantial theoretical interest.
\begin{quote}
\textbf{Open problem.} Why is standard time frequency analysis applied to $\left|f(t)\right|$ better able to recover the fundamental underlying frequency? Why does the `rectification' trick work?
\end{quote}

Abstracting the trick, we may interpret the rectification trick as the application of a   {\em nonlinear activation function} $g(t) = |t|$: instead of analyzing $f(t)$, one analyzes $g(f(t))$. Naturally, there are many other choices of the activation function \cite{Rasamoelina2020}, like the commonly used Rectified Linear Unit (ReLU) function $g(t) = \mbox{ReLU}(f(t))$ and other widely applied activation functions from the theory of neural networks that are at our disposal. In practice, ReLU seems to work and work roughly as well as the absolute value. We do not know of any theoretical support for any of these functions and can now formulate the main question that motivates our paper.
\begin{quote}
\textbf{Open problem.} Which nonlinear activation function $g:\mathbb{R} \rightarrow \mathbb{R}$ leads to the `best' recovery of the fundamental frequency?
\end{quote}

One of the main points of our paper is to discuss the problem from the point of view of classical Fourier analysis and to propose, based on that, a somewhat unorthodox choice, provide theoretical support and discuss its behavior in practice.

\section{Main Result}
\subsection{An adaptive activation function.}
We define, for $\varepsilon > 0$ small, an adaptive activation function. Introducing, for $0 < \varepsilon < 1$, the function $h_{\varepsilon}: [-1,1] \rightarrow \mathbb{R}_{\geq 0}$
$$ h_{\varepsilon}(x) =  \frac{1}{1 - (1 - \varepsilon)|x|},$$
we propose to run standard time frequency on the normalized signal
$$
h_{\varepsilon}\left( \frac{f(t)}{\|f\|_{L^{\infty}}}\right) = \frac{1}{1 - (1 - \varepsilon)\frac{|f(t)|}{\|f\|_{L^{\infty}}}} 
$$
and expect it to enhance the fundamental component of $f$. This activation function is adaptive to the input signal since it depends on the $L^\infty$ norm of the input signal. Note that $h_{\varepsilon}$ is also $2\pi-$periodic. So by construction, it inherits the periodicity of $f$. The question is now: can we determine the strength of the fundamental component of $h_{\varepsilon}$?
To provide a theoretical analysis this question, we put some assumptions.
\begin{enumerate}
\item We assume that the function $f:[0,2\pi] \rightarrow \mathbb{C}$ is of the type
$$ 0 \neq f(t) = \sum_{k=1}^{n} a_k e^{i m_k t} \quad a_k \in \mathbb{C}, m_k \in \mathbb{N},~ 0 < m_1 < m_2 < \dots < m_n.$$
In particular, this function will have fundamental frequency $\gcd(m_1, \dots, m_n)$.
\item We assume that the function $g(t) = \left|f(t)\right|$ assumes its global maximum in finitely many points $\left\{ t_1, \dots, t_m \right\} \subset [0,2\pi]$ and
\item that these maxima are non-degenerate: $g''(t_i) < 0$ for all $1 \leq i \leq m$.
\end{enumerate}

We can now state our result which shows that, as $\varepsilon \rightarrow 0$, the method will indeed recover the fundamental frequency (which we 
expect to be 1 in the generic setting).

\begin{theorem}
We have, as $\varepsilon \rightarrow 0$,
$$\int_0^{2\pi} h_{\varepsilon}\left( \frac{f(t)}{\|f\|_{L^{\infty}}}\right) e^{it} dt =   \frac{1}{\sqrt{\varepsilon}} \sum_{j=1}^m \frac{ \pi \cdot e^{i t_j}}{\sqrt{ - \frac{g''(t_j)}{2 \|g\|_{L^{\infty}}}}} +  \mathcal{O}\left(\frac{1}{\varepsilon^{1/4}}\right).$$
\end{theorem}

\textit{Remarks.}
\begin{enumerate}
\item We note that for `generic' signals (say, random signals picked from a suitable probability distribution), we expect that there is a unique maximum $t_1$: in that case, the leading order expansion is always guaranteed to be of order 
$$ \widehat{h_{\varepsilon}}(1) =   \frac{1}{\sqrt{\varepsilon}}  \frac{\pi \cdot e^{i t_1}}{\sqrt{ - \frac{g''(t_1)}{2 \|g\|_{L^{\infty}}}}} + \mathcal{O}\left(\frac{1}{\varepsilon^{1/4}}\right)$$
and we are guaranteed that the fundamental frequency is recovered.
\item If there is more than one maximum, then we only fail at recovering the fundamental frequency if 
$$ \sum_{j=1}^m \frac{e^{i t_j}}{\sqrt{ - g''(t_j)}}  = 0.$$
This is something that we do not expect to happen for generic signal (it corresponds to a precise algebraic identity). 
\item If the frequencies share a nontrivial greatest common divisor 
$$G=\gcd(m_1, \dots, m_n) > 1,$$ then the maxima of $g(t)$ will be $G-$periodic on the unit circle and the sum will vanish. The Theorem can then be applied to the function $g(t) = f(t/G)$ to conclude that the fundamental frequency $G$ will be recovered.
\end{enumerate}

\subsection{General Activation Functions} The purpose of this section is to give some perspective on the problem from the point of view of classical Fourier analysis. We note that the intrinsically nonlinear nature of the problem does pose an interesting challenge, nonetheless, there are some natural considerations partially inspired by a classical paper of Rudin \cite{rudin} that seem like they might be relevant.
We recall that our goal is, given a function $f$ of the type
$$ 0 \neq f(t) = \sum_{k=1}^{n} a_k e^{i m_k t} \qquad a_k \in \mathbb{C}, m_k \in \mathbb{N}, 0 < m_1 < m_2 < \dots < m_n$$
to come with a nonlinear function $g: \mathbb{R} \rightarrow \mathbb{R}$ such that $g(f(t))$ has a clear frequency contribution at
$ \gcd\left(m_1, \dots, m_n\right)$. The main difficulty in practice is, of course, that we do not actually know the precise formula of $f$ when only $f$ is given, and have no
real idea what the actual frequencies $m_1, \dots, m_n$ might be. 
We simplify things by assuming that the nonlinearity $g$ can be expanded into a Taylor series
$$ g(x) = b_0 + b_1 x + b_2 x^2 + b_3 x^3 + \dots.$$
Then each individual power applied to $f$ can be properly analyzed: note that
$$ \left(  \sum_{k=1}^{n} a_k e^{i m_k t} \right)^k = \sum_{\ell \in \mathbb{Z}} \left(\sum_{m_{i_1} + \dots + m_{i_r} = \ell} a_{i_1} a_{i_2} \dots a_{i_r} \right) e^{i \ell t}.$$
This shows that large powers naturally lead to a function whose frequencies is comprised of sums of the individual frequencies. If we replace
powers by powers of the absolute value, then at least for even powers the identity $|x|^{2k} = x^k \overline{x}^k$ leads to a way of writing
$$ X = \left|  \sum_{k=1}^{n} a_k e^{i m_k t} \right|^{2k}$$
as
$$X = \sum_{\ell_1, \ell_2 \in \mathbb{Z}} \left(\sum_{\scriptsize m_{i_1} + \dots + m_{i_r} = \ell_1} a_{i_1} \dots a_{i_r} \right)
\left(\sum_{m_{i_1} + \dots + m_{i_r} = \ell_2} \overline{ a_{i_1}}   \dots \overline{a_{i_r}} \right)  e^{i (\ell_1-\ell_2) t}.$$
There is a nontrivial contribution coming from the coefficients $(a_k)_{k \in \mathbb{Z}}$, however, if we only consider the support of the frequencies, then a rather
clear picture emerges: denoting the set of frequencies that make up $f$ by
$$ M = \left\{m_1, \dots, m_n\right\},$$
we expect that
$$ \mbox{supp}~ \mathcal{F} |f|^{2k} = \underbrace{(M + M + \dots +M)}_{k~\mbox{\tiny times}} -  \underbrace{(M + M + \dots +M)}_{k~\mbox{\tiny times}}.$$
However, we also expect, for $k$ sufficiently large, that the numbers that can be appropriately written as sums and differences of elements in $M_n$ have a nice limiting behavior and (for a suitabye lenient definition of $\lim$)
$$ \lim_{k \rightarrow \infty} \quad \underbrace{(M + M + \dots +M)}_{k~\mbox{\tiny times}} -  \underbrace{(M + M + \dots +M)}_{k~\mbox{\tiny times}}  = \gcd(m_1, \dots, m_n) \cdot \mathbb{Z}.$$
These heuristic considerations suggest the following rough guidelines:
\begin{enumerate}
\item The Fourier transform of $g(f(t))$ depends on the Taylor expansion of $g$.
\item If the function $g$ is smooth, then the Taylor expansion will have rapidly decaying coefficients and sum sets of the type 
$$  \underbrace{(M + M + \dots +M)}_{k~\mbox{\tiny times}} \quad \mbox{only contribute for small}~k.$$
A choice like $g(x) = \cos{(x)}$ will not lead to a good detection function.
\item This suggests functions $g$ for which the Taylor expansion decays slowly:
functions that are discontinuous or have discontinuous first derivative.
\item Moreover, introducing the absolute value $g(|x|)$ leads to sum-difference sets which might generally lead to better results.
\end{enumerate}

These principles naturally suggest functions like $g(x) = |x|$ or $g(x) = \mbox{ReLu}(x)$. However, we note that both of these functions
are actually continuous which leads to a decay of coefficients in the Taylor expansion. Indeed, considering all these principles, we should
focus our attention on functions of the form
$$ g(x) = b_0 + b_1 |x| + b_2 |x|^2 + b_3 |x|^3 + \dots$$
for which the $(b_j)_{j=0}^{\infty}$ decay slowly. A particularly natural candidate is
$$ \frac{1}{1-|x|} = 1 + |x| + |x|^2 + \dots$$
and this is how we arrive at our proposed construction. As we will show in Section 3, this construction is not only well motivated by the aforementioned
arguments but also naturally incorporates a second approach which is related to asymptotic analysis and encapsulated by Theorem 1. 

\section{Proof of the Theorem} 
\begin{proof}
We now give a proof of the result.
 It follows from the assumptions that each local maximum of $g$ behaves locally around the maximum like a parabola (since $g''$ does not vanish and is negative). 
This allows us to treat all the points where $f$ assumes a maximum in isolation. We start by considering a point $t_0$ such that $g(t_0) \neq \|g\|_{L^{\infty}}$. In that case, we observe that 
$$h_{\varepsilon}(t_0) = \frac{1}{1 - \frac{1 - \varepsilon}{\|f\|_{L^{\infty}}} |f(t_0)|} \leq \frac{1}{1- \frac{|f(t_0)|}{\|f\|_{L^{\infty}}}} \leq \frac{\|f\|_{L^{\infty}}}{\|f\|_{L^{\infty}} - |f(t_0)|}.$$
This quantity is uniformly bounded as $\varepsilon \rightarrow 0$. Since our main result is about the asymptotic growth of an integral over a bounded region as $\varepsilon \rightarrow 0$ and since we are only interested in terms at scale $\sim \varepsilon^{-1/2}$ or larger, we may disregard points where $h_{\varepsilon}(t)$ remains bounded as $\varepsilon \rightarrow 0$. We can now proceed as follows: since there are only finitely many points $t_1, t_2, \dots, t_m$ that achieve the global maxima of $g$ and since all of them are non-degenerate, there exists $\delta > 0$ such that the preimage
$$g^{-1}\left[\|f\|_{L^{\infty}} - \delta, \|f\|_{L^{\infty}}\right] = \bigcup_{j=1}^{m} I_j$$
 can be written as the union of $m$ disjoint intervals. By making $\delta$ sufficiently small, we can also ensure that each of the intervals contains a 
 point in which $|f(t)|$ assumes its global maximum.
 Moreover, we can also infer that, as $\delta \rightarrow 0$, these intervals scale asymptotically like
 $ |I_j| \sim \delta^{-1/2}$, where the implicit constant depends on the value of the second derivative in $t_j$. 

\begin{center}
\begin{figure}[h!]
\begin{tikzpicture}
\draw [dashed] (-0.5,1) -- (8,1);
\draw [dashed] (-0.5,0) -- (8,0);
\node at (8.7, 1) {$\|f\|_{L^{\infty}}$};
\node at (9, 0) {$\|f\|_{L^{\infty}} - \delta$};
\draw [thick] (0,0) to[out=80, in=180] (1,1) to[out=0, in=110] (2,0);
\draw [thick] (3,0) to[out=40, in=180] (4,1) to[out=0, in=140] (5,0);
\draw [thick] (6,0) to[out=80, in=180] (6.8,1) to[out=0, in=110] (7.6,0);
\node at (1, -0.3) {$I_1$};
\node at (4, -0.3) {$I_2$};
\node at (6.8, -0.3) {$I_3$};
\end{tikzpicture}
\caption{A sketch of the decomposition.}
\end{figure}
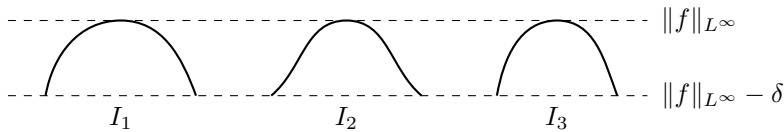
\end{center}

Then, uniformly as $\varepsilon \rightarrow 0$, we have
\begin{align*}
 &\int_0^{2\pi} h_{\varepsilon}(t) 1_{\left\{|f| < \|f\|_{L^{\infty}} - \delta\right\}}(t) dt \leq \int_0^{2\pi} \frac{\|f\|_{L^\infty}}{\|f\|_{L^\infty}-g(t)} 1_{\left\{|f| < \|f\|_{L^{\infty}} - \delta\right\}}(t) dt\\
 \leq&\,\int_0^{2\pi} \frac{\|f\|_{L^\infty}}{\delta} 1_{\left\{|f| < \|f\|_{L^{\infty}} - \delta\right\}}(t) dt\leq  \frac{2 \pi \|f\|_{L^{\infty}}}{\delta}\,. 
 \end{align*} 
Note that $ 2 \pi \|f\|_{L^{\infty}}/\delta$ is independent of $\epsilon$.
This shows that contributions at scale $\varepsilon^{-1/2}$ can only come from $I_1, I_2, \dots, I_m$. Let us now assume that $g$ assumes a
global maximum in $t_1 \in I_1$.
Then, locally around $t_1$, we have the expansion
$$ g(t) = \|g\|_{L^{\infty}} + \frac{g''(t_1)}{2} (t-t_1)^2 + \mathcal{O}((t-t_1)^3).$$
Therefore, we have, locally around $t_1$,
\begin{align}
 h_{\varepsilon}(t) &= \frac{1}{1 - (1-\varepsilon) \left( 1 + \frac{g''(t_1)}{2 \|g\|_{L^{\infty}}} (t-t_1)^2 + \mathcal{O}((t-t_1)^3)\right)}\nonumber \\
 &= \frac{1}{\varepsilon - (1-\varepsilon)  \frac{g''(t_1)}{2 \|g\|_{L^{\infty}}} (t-t_1)^2 + \mathcal{O}((t-t_1)^3)}.\label{expansion of heps near t1}
  \end{align}
  Note that $g'' < 0$. Thus, the denominator is, up to second order, growing away from $t_1$.
Our next ingredient will be the identity
\begin{align}
\forall~A,B,C > 0 \qquad \quad \int_{|t| > C} \frac{1}{A + B t^2} dt = \frac{\pi - 2\arctan\left(\frac{\sqrt{B}}{\sqrt{A}} C\right)}{\sqrt{A B}}\label{key integration part1}
\end{align}
In particular, letting $C \rightarrow 0$, we recover
$$\forall~A,B > 0 \qquad \quad \int_{\mathbb{R}} \frac{1}{A + B t^2} dt = \frac{\pi}{\sqrt{A B}}.$$
In our case, we have $A = \varepsilon>0$, while 
$$ B = -(1-\varepsilon)  \frac{g''(t_1)}{2 \|g\|_{L^{\infty}}}>0$$
converges to a fixed constant as $\varepsilon$ gets small (and is, in particular, bounded away from 0). This implies
\begin{align*}
\int_{-\varepsilon^{1/4}}^{\varepsilon^{1/4}} \frac{1}{\varepsilon + B t^2} dt &=\int_{-\infty}^\infty \frac{1}{\varepsilon + B t^2} dt-\int_{|t|\geq \varepsilon^{1/4}} \frac{1}{\varepsilon + B t^2} dt \\
&=  \frac{\pi}{\sqrt{\varepsilon B}} + \mathcal{O}_{}(B^{-1/2}\varepsilon^{-1/4})
\end{align*}
when $\varepsilon$ is sufficiently small. The second equality comes from \eqref{key integration part1}. Indeed, by choosing $C=\varepsilon^{1/4}$, we see that the term $\sqrt{B}/\sqrt{A}C$ is of order $\varepsilon^{-1/4}$, and hence 
$$\pi-2\arctan\left(\frac{\sqrt{B}}{\sqrt{A}} C\right) \lesssim \frac{\sqrt{A}}{\sqrt{B}}C^{-1}$$ whenever $\varepsilon$ is sufficiently small. Combined with the numerator, we obtain the desired control of $\int_{|t|\geq \varepsilon^{1/4}} 1/(\varepsilon + B t^2) dt $.
From this control, we deduce 
$$ \int_{I_1} h_{\varepsilon}(t) dt =  \frac{\pi}{\sqrt{\varepsilon}} \frac{1}{\sqrt{ -\frac{g''(t_1)}{2 \|g\|_{L^{\infty}}}}} +\mathcal{O}\left(\frac{1}{\varepsilon^{1/4}}\right)$$
when $\varepsilon$ is sufficiently small. Indeed, since $B \sim 1$, the error, $g''(t_1)/(2 \|g\|_{L^{\infty}}) + \mathcal{O}(t-t_1)$ in \eqref{expansion of heps near t1}, is controlled when $\varepsilon$ is sufficiently small.
Using continuity of $e^{it}$ and summing over the $m$ intervals, we arrive our conclusion
$$ \int_{0}^{2\pi} h_{\varepsilon}(t) e^{i t} dt = \frac{\pi}{\sqrt{\varepsilon}} \sum_{j=1}^m \frac{e^{i t_j}}{\sqrt{ - \frac{g''(t_j)}{2 \|g\|_{L^{\infty}}}}} + \mathcal{O}\left(\frac{1}{\varepsilon^{1/4}}\right).$$
\end{proof}

\section{Numerical results}

All Matlab codes are available in \url{} for the reproducibility purposes. 

\subsection{Synthetic data.}

\begin{figure}[htb!]\centering
	\includegraphics[trim=0 0 0 0, clip, width=0.99\textwidth]{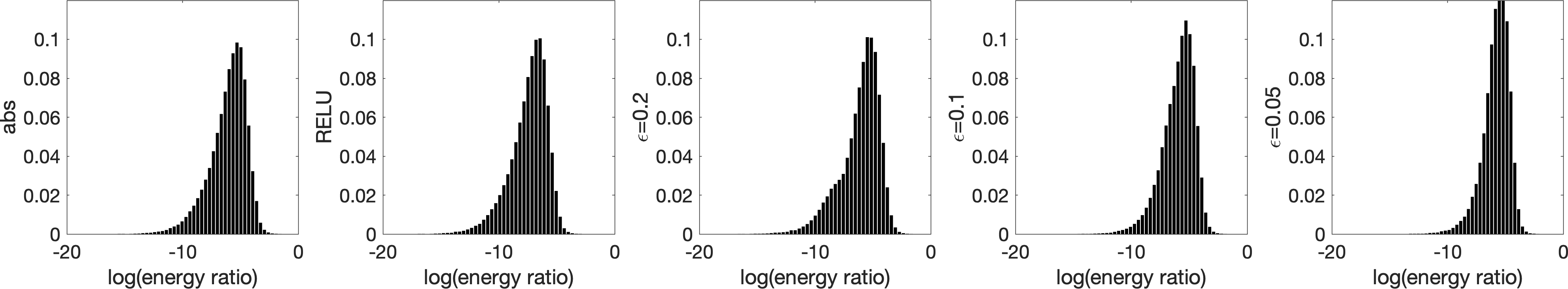}
	\caption{Effect of different activation functions on fundamental component enhancement. 
	The histograms of the fundamental component enhancement by 5 different activation functions, including $|\cdot|$, RELU, $h_{0.2}$, $h_{0.1}$ and $h_{0.05}$, are shown from left to right. \label{Figure6synthetic}}
\end{figure}

We generate synthetic data to investigate how different activation functions enhance the fundamental component. Fix the sampling rate to be $512$ Hz. First, randomly select an integer $K$ between $5$ and $100$ according to a uniform distribution. Then, randomly select $K$ integers between $2$ and $250$ following the probability density function $Ce^{-(l/100)^2}$, where $C$ is the normalization constant, so that the chosen integers have 1 as the greatest common divisor. Denote the selected integers as $\{j_k\}_{k=1}^K$. Then, set the function $f(t)$ as 
\[
f(t)=\sum_{k=1}^K a_k\cos(2\pi j_k t+\phi_k)\,,
\]
where $a_k$, $k=1,\ldots,K$, are identically and independently (i.i.d.) sampled from $(0,1]$, and $\phi_k$, $k=1,\ldots,K$, are i.i.d. sampled from $(0,2\pi]$. Note that by construction, $f(t)$ is a toy example as that shown in \eqref{Simple model harmonic 1} with $\hat{f}(1)=0$ but $1$ as the fundamental frequency. For a given activation function $g(t)$, we evaluate the resulting  fundamental component enhancement by evaluating the {\em fundamental component energy ratio}, defined as
\begin{equation}\label{def ER FC}
\frac{|\widehat{g\circ f}(1)|^2}{\sum_{l=1}^{256}|\widehat{g\circ f}(l)|^2}\,.
\end{equation}
Then, repeat the above procedure for $10^5$ times.
The results of different activation functions are shown in Figure \ref{Figure6synthetic} as histograms of fundamental component energy ratios.
Quantitatively, the median and median absolute deviation of the fundamental component energy ratio over $10^5$ times are
0.28\% and 0.46\% for the rectification,
0.07\% and 0.12\% for the RELU,
0.29\% and 0.46\% for $h_{0.2}$,
0.31\% and 0.42\% for $h_{0.1}$, and
0.33\% and 0.33\% for $h_{0.05}$.

\subsection{Semi-real Medical Example.}

\begin{figure}[htb!]\centering
	\includegraphics[trim=0 0 0 0, clip, width=0.99\textwidth]{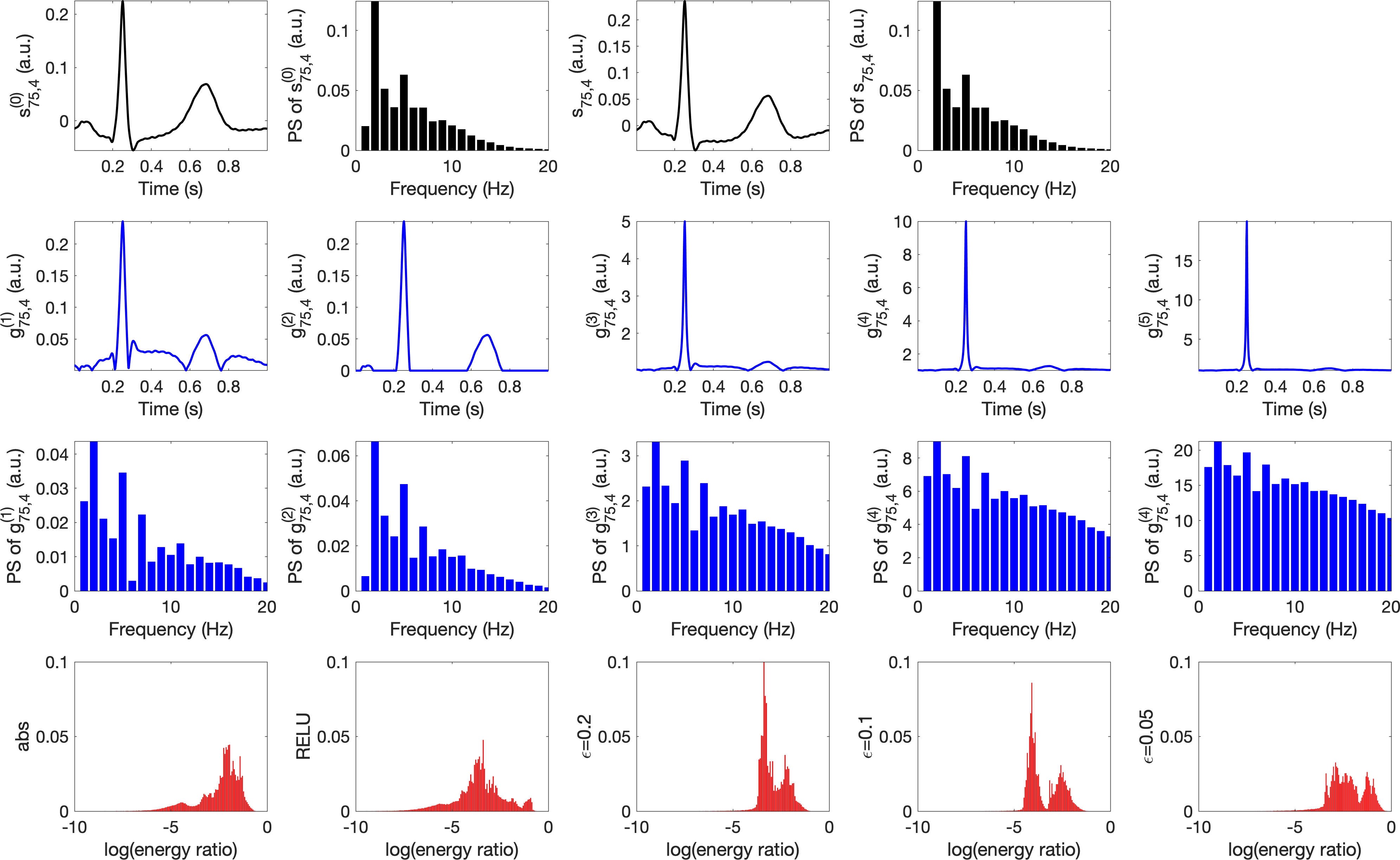}
	\caption{Illustration of the effect of different activation functions on the fundamental component enhancement. 
	The first row shows one cardiac cycle (left two subplots) and the simulated cardiac cycle without the fundamental component (right two subplots) from the CinC2013 database.
	The second row shows the resulting signal after applying
	5 different activation functions, including $|\cdot|$, RELU, $h_{0.2}$, $h_{0.1}$ and $h_{0.05}$ from left to right.
	The third row shows the associated power spectra (PS) of those 1-periodic functions shown on the second row. 
	The fourth row shows the histograms of energy ratio of the fundamental component after applying different activation functions, where the x-axis is the log of the energy ratio. The (medians, median absolute deviation) of the energy ratios are (11.7\%, 7\%), (2.9\%, 5.2\%), (4.4\%, 3.9\%), (2.2\%, 3.3\%) and (9.8\%, 10.8\%) respectively.  a.u. means arbitrary unit. \label{Figure7simu}}
\end{figure}

\begin{figure}[htb!]\centering
	\includegraphics[trim=0 0 0 0, clip, width=0.998\textwidth]{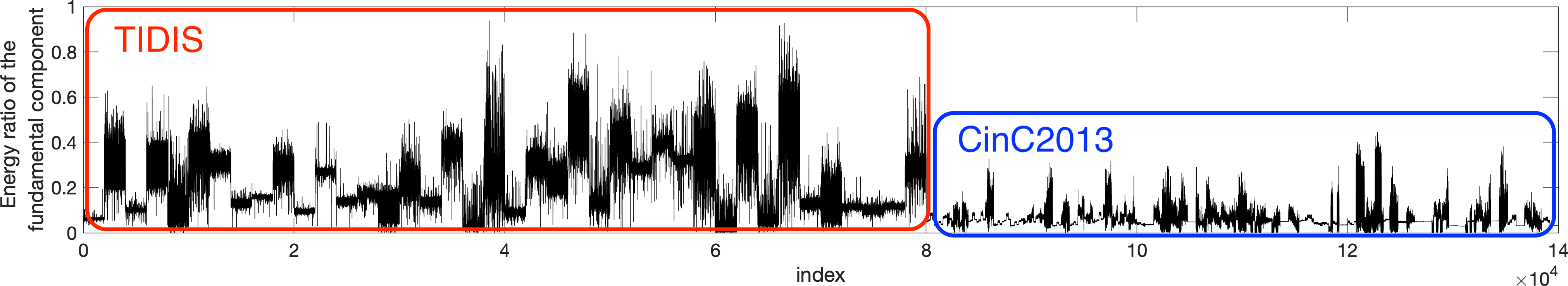}
	\includegraphics[trim=0 0 0 0, clip, width=0.998\textwidth]{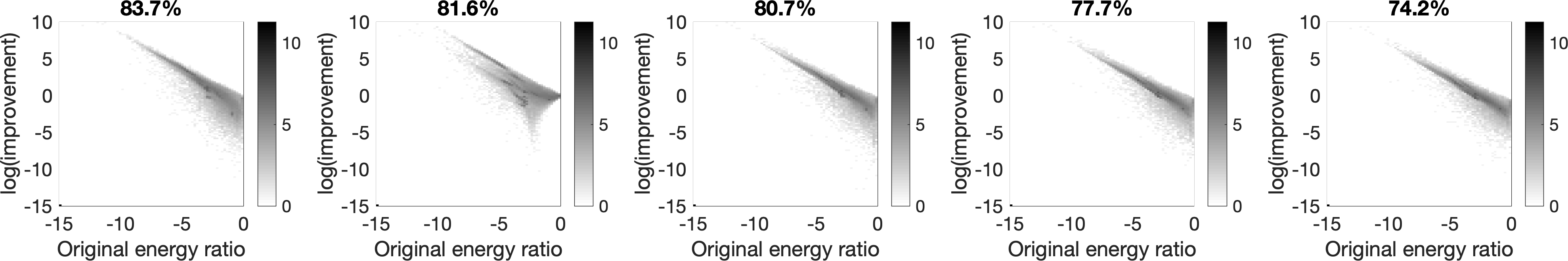}
	\caption{Top: energy ratio of the fundamental component, $r^{(0)}_{ik}$ of all cardiac cycles. The first part marked by the red box comes from the TIDIS database, and the second part marked by the blue box comes from the CinC2013 database. Bottom: two-dimensional histograms of energy ratio of the fundamental component after applying different activation functions, $|\cdot|$, RELU, $h_{0.2}$, $h_{0.1}$ and $h_{0.05}$ from left to right, where the x-axis is $\log(r^{(0)}_{ik})$, and the y-axis is $\log(r^{(j)}_{ik})/\log(r^{(0)}_{ik}))$. The portion of cardiac cycles that are enhanced by the nonlinear activation function is listed in the title. \label{Figure8simu}}
\end{figure}

We generate a semi-real dataset and compare the proposed activation function with $\epsilon=0.2$, $\epsilon=0.1$ and $\epsilon=0.05$ with the usually applied rectification and RELU activation function. 
We consider two databases. The first one is the noninvasive trans-abdominal ECG database from the PhysioNet Computing in Cardiology Challenge 2013 (CinC2013)\footnote{\url{https://physionet.org/content/challenge-2013/1.0.0/}}. There are in total 75 recording, each recording contains 4 channels, and each record lasts for 1 minute. The second one is the Taiwan Integrated Database for Intelligent Sleep (TIDIS)\footnote{\url{https://tidis.org}}. There are in total 20 whole night polysomnogram recordings, each recording contains 2 channels, and each record lasts for about 6 hours. 
For each recording in CinC2013, we detect the maternal R peaks by applying the standard R peak detection algorithm to the mean of 4 channels; for each recording in TIDIS, we apply the same standard R peak detection algorithm to the mean of two channels. For the $k$-th channel in the $i$-th recording, we randomly generate an 1-periodic function, denoted as $s_{ik}^{(0)}(t)$, by averaging $10$ randomly selected cardiac cycles. 
Then we intentionally remove the fundamental component via the Fourier transform, and denote the resulting 1-periodic function as $s_{ik}(t)$ that is sampled at 512 Hz. Next, we apply the $j$-th activation function to $s_{ik}(t)$, and denote the resulting signal as $g^{(j)}_{ik}(t)$. To evaluate how much the fundamental is activated, we record the energy ratio of the fundamental component defined in \eqref{def ER FC}, denoted as $r^{(j)}_{ik}$.
For the $k$-th channel in the $i$-th recording, the above procedure is repeated $100$ times. The histograms of all collected energy ratios for different activation functions are shown in Figure \ref{Figure7simu}. We see that the proposed activation function gives a delta-like signal with the peak centered at the maximal value. Hence, the power spectrum is flatter compared with both rectification and RELU. This explains why the median of the energy ratio is smaller if we apply the proposed activation function since the proposed activation function tends to flatten the spectrum, particularly when $\epsilon$ is small.
We shall mention that the main difference between these two databases is that the cardiac cycle is recorded from the abdomen in the CinC2013 database, and from the chest in the TIDIS database. Thus the cardiac cycles are of different morphology, and hence different energy ratio of the fundamental component. Note that after the averaging, the fetal ECG impact is reduced, and we obtain a reasonably clean maternal ECG cycles. See Figure \ref{Figure8simu} for more details.

To have a complete picture of the fundamental enhancement, we repeat the same data preparation procedure above, except removing the fundamental component; that is, we simulate the practical situation that the fundamental component may or may not exist, and its strength may or may not be weak. In addition to record the energy ratio, also denoted as $r^{(j)}_{ik}$ by a bit abuse of notation, after applying the activation function, we also record the energy ratio before, denoted as $r^{(0)}_{ik}$. In Figure \ref{Figure6simu},  we plot the 2-dim histogram of the distribution of $\{(\log(r^{(0)}_{ik}), \log(r^{(j)}_{ik})/\log(r^{(0)}_{ik}))\}$. By definition, when $\log(r^{(j)}_{ik})/\log(r^{(0)}_{ik}))>1$, we obtain an enhancement of the energy ratio of the fundamental component. The portion of cardiac cycles that the fundamental component strength is enhanced is recorded in the title. It is clear that all activation functions enhance the energy ratio of the fundamental component when it is weak.

Note that at the first glance, our proposed activation function seems to provide a worse result due to the smaller portion of enhanced cardiac cycles. However, we should emphasize that due to the ``flattening'' nature of our activation function, the energy ratio of the fundamental component is smaller. 

\subsection{Medical Example.}

\begin{figure}[htb!]\centering
	\includegraphics[trim=0 0 0 0, clip, width=0.9\textwidth]{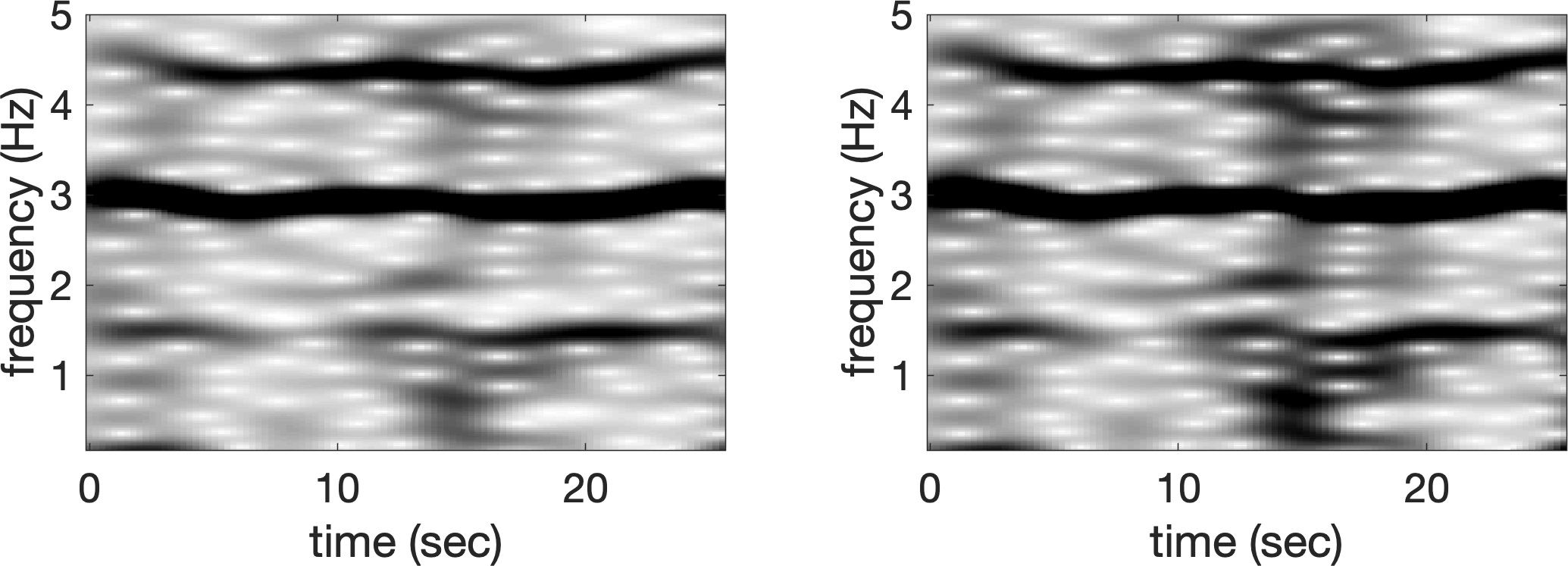}
	\caption{The spectrogram of the PCG signal composed with RELU (left) and $h_{0.01}$, where we use the same window and the same dynamic range as those shown in Figure \ref{fig:pcg}(b). \label{Figure9PCG}}
\end{figure}

\begin{figure}[htb!]\centering
	\includegraphics[trim=0 0 0 0, clip, width=0.98\textwidth]{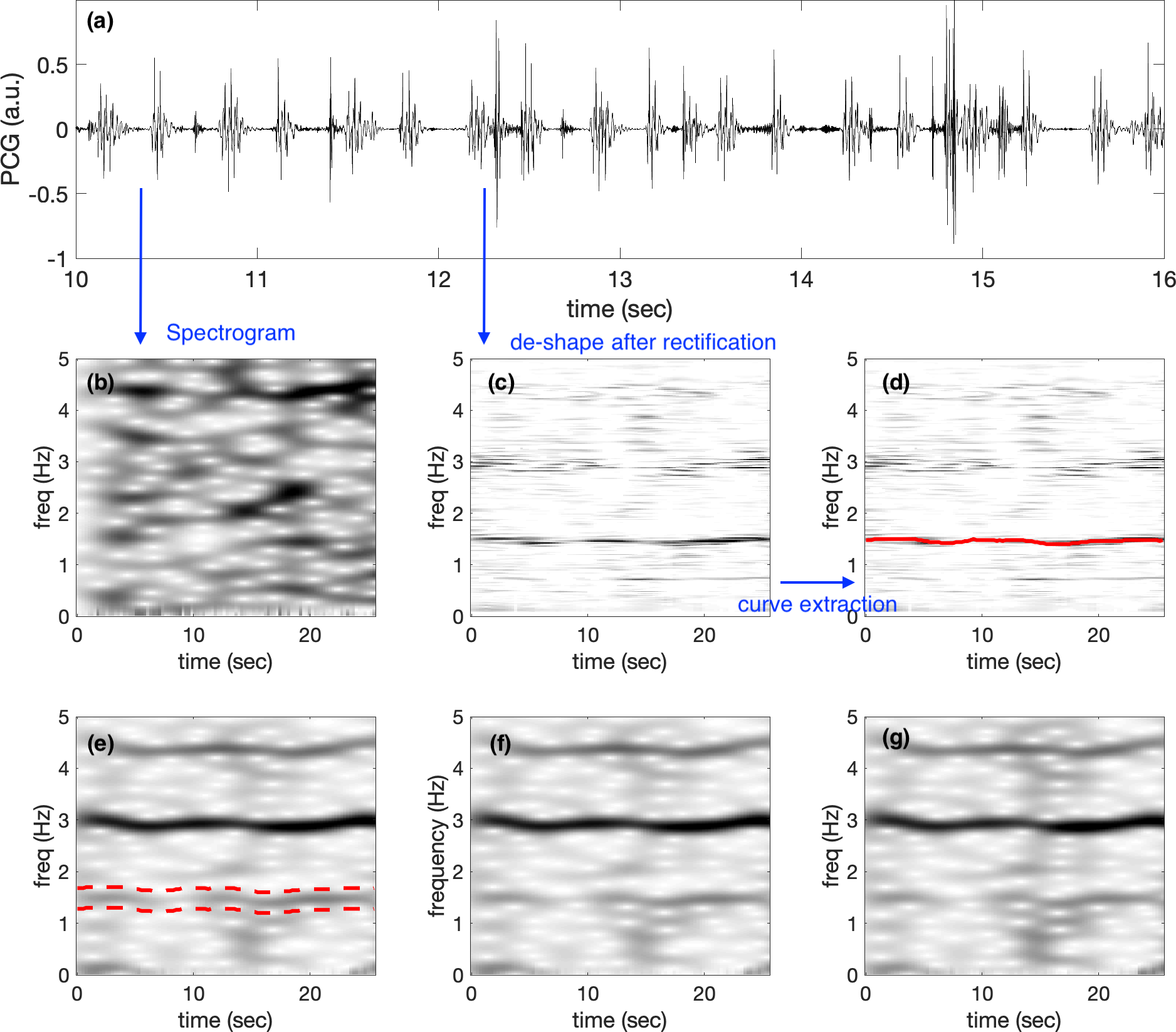}
	\caption{(a) a PCG signal $f(t)$. (b) the spectrogram of $f(t)$, where we show only the frequency range from 0 Hz to 5 Hz. (c) the de-shape spectrogram of $|f(t)|$, where we show only the frequency range from 0 Hz to 5 Hz. (d) the de-shape spectrogram superimposed with a red curve representing the detected instantaneous frequency.  (e)-(g) are spectrograms associated with $|f(t)|$, $\texttt{RELU}(f(t))$ and $h_{0.1}$. In (e), two red curves are superimposed representing the spectral band that we use to calculate the energy ratio of the fundamental component associated with the detected instantaneous frequency. In all time-frequency representations, the dynamic range is set to be the 0\% and 99.95\% percentiles of all entries. \label{Figure10PCG}}
\end{figure}

Next, we consider the PCG signal as a real example. In Figure \ref{Figure9PCG}, we show the spectrogram of the PCG signal shown in shown in Figure \ref{Figure3PCG}(a) composed with RELU (left) and $h_{0.1}$. We could see that both RELU and $h_{0.1}$ provide the fundamental component information; that is, there is an identifiable curve around 2Hz. To have a quantification of this finding, we analyze the maternal PCG signals in the Shiraz University Fetal Heart Sounds Database\footnote{\url{https://physionet.org/content/sufhsdb/1.0.1/}}. There are in total 92 recordings. 
See Figure \ref{Figure10PCG}(a) for a PCG signal different from that in Figure \ref{fig:pcg}(a). It is clear that the spectrogram provides limited information about the heart rate. See Figure \ref{Figure10PCG}(b) and recall Figure \ref{fig:pcg}(b) for an illustration. To quantify the fundamental component enhancement performance, we carry out the following steps. First, run the de-shape algorithm \cite{LinLiWu2018} to determine the instantaneous frequency (IF) of the PCG. See Figure \ref{Figure10PCG}(c-d) for the de-shape spectrogram and the detected IF. Denote the estimated IF as $\phi'(t)$. The detected heart rate is confirmed by reading and hearing the signal. Second, determine the energy ratio around the band $[\phi'(t)-0.2, \phi'(t)+0.2]$ Hz; that is, 
\[
R:=\frac{\int_0^T \int_{\phi'(t)-0.2}^{\phi'(t)+0.2} |V(t,\xi)|^2d\xi }{\int_0^T \int_{1/T}^{U} |V(t,\xi)|^2d\xi}\,,
\]
where $|V|^2$ is the spectrogram of the signal of length $T$ with the sampling rate $1/U$. Here, since our focus is the heart rate, we downsample the signal to $U=100$ Hz. 
Over the 92 recordings, the mean and standard deviation of $R$ determined from the original signal, the rectified signal, the RELU activated signal and $h_{0.1}$ are $0.03\%\pm 0.01\%$, $2.33\%\pm 0.56\%$, $1.64\%\pm 0.41\%$, and $1.25\%\pm 0.28\%$. See Figure \ref{Figure10PCG} for an example of the overall procedure.

\end{document}